\input amstex
\documentstyle{amsppt}
\document
\topmatter
\title
Holomorphically pseudosymmetric K\"ahler metrics on $\Bbb{CP}^n$.
\endtitle
\author
W{\l}odzimierz Jelonek
\endauthor

\abstract{The aim of this paper is to present   examples of
K\"ahler holomorphically pseudosymmetric K\"ahler metrics on the
complex projective spaces $\Bbb{CP}^n$. }
\thanks{MS Classification: 53C55,53C25. Key words and phrases:
K\"ahler manifold, holomorphically pseudosymmetric K\"ahler
manifold, QCH K\"ahler manifold}\endthanks
 \endabstract
\endtopmatter

\define\r{\rightarrow}

\define\n{\nabla}
\define\om{\omega}

\define\k{\diamondsuit}
\define\th{\theta}

\define\a{\alpha}

\define\lb{\lambda}

\define\1{D_{\lb}}
\define\2{D_{\mu}}
\define\0{\Omega}

\define\De{\Cal D}

\define\m{(M,g,J)}

{\bf 1. Introduction.}  The aim of this paper is to give new
examples of  holomorphically pseudosymmetric K\"ahler metrics on
complex projective spaces.  The holomorphically pseudosymmetric
K\"ahler manifolds were defined by Z. Olszak in [O-1] in 1989 and
studied in [D], [H], [Y]. A K\"ahler manifold $\m$ is called
holomorphically pseudosymmetric (HP K\"ahler) if its curvature
tensor $R$ satisfies the condition
$$R.R=\phi\Pi.R,$$
where $R,\Pi$ act as  derivations of the tensor algebra,
$\Pi(U,V)X=\frac14(g(V,X)U-g(U,X)V+g(JV,X)JU-g(JU,X)JV-2g(JU,V)JX)$
is the K\"ahler type curvature tensor of constant holomorphic
sectional curvature and $\phi\in C^{\infty}(M)$ is a smooth
function. A Riemannian manifold $(M,g)$ is called semisymmetric if
$R.R=0$ (see [Sz-1], [Sz-2]). Until recently  examples of compact,
holomorphically pseudosymmetric and not semisymmetric, K\"ahler
manifolds were not known. The first compact, simply connected
examples of HP K\"ahler manifolds which are not semisymmetric were
given by the author in [J-2]. These manifolds turned out to be QCH
K\"ahler manifoldas. The QCH K\"ahler manifolds are the K\"ahler
manifolds admitting a smooth, two-dimensional, $J$-invariant
distribution $\De$ whose holomorphic curvature
$K(\pi)=R(X,JX,JX,X)$ of any $J$-invariant $2$-plane $\pi\subset
T_xM$, where $X\in \pi$ and $g(X,X)=1$, depends only on the point
$x$ and the number $|X_{\De}|=\sqrt{g(X_{\De},X_{\De})}$, where
$X_{\De}$ is the orthogonal projection of $X$ on $\De$. In this
case  we have
$$R(X,JX,JX,X)=\phi(x,|X_{\De}|)$$ where $\phi(x,t)=a(x)+b(x)t^2+c(x)t^4$ and
 $a,b,c$ are smooth functions on $M$. Also $R=a\Pi+b\Phi+c\Psi$
 for certain curvature tensors $\Pi,\Phi,\Psi\in \bigotimes^4\frak X^*(M)$
  of K\"ahler type (see [G-M-1],[G-M-2], [J-1]). In the present
  paper we construct QCH K\"ahler metrics on the manifold $\Bbb{CP}^n-\{p_0\}$,
  where $p_0\in \Bbb{CP}^n$ and show that these metrics extend smoothly
    to HP K\"ahler metrics on the complex projective
  spaces  $\Bbb{CP}^n$.
  \medskip
{\bf 2.  QCH K\"ahler manifolds.}  By $h$ we shall denote the
tensor $h=g\circ (p_{\De}\times p_{\De})$, where $p_{\De}$ are the
orthogonal projections on $\De$.  By $\0=g(J\cdot,\cdot)$ we shall
denote the K\"ahler form of $\m$, by $\om$  the K\"ahler form of
$\De$ i.e. $\om(X,Y)=h(JX,Y)$. We shall recall some results proved
by Ganchev and Mihova in [G-M-1]. Let
$R(X,Y)Z=([\n_X,\n_Y]-\n_{[X,Y]})Z$ and let us write
$$R(X,Y,Z,W)=g(R(X,Y)Z,W).$$ We shall identify $(1,3)$ tensors
with $(0,4)$ tensors in this way. If  $R$ is the curvature tensor
of a QCH K\"ahler manifold $\m$, then
$$R=a\Pi+b\Phi+c\Psi,$$
where  $a,b,c\in C^{\infty}(M)$ and    $\Pi$ is the standard
K\"ahler tensor of constant holomorphic curvature
$$\gather \Pi(X,Y,Z,U)=\frac14(g(Y,Z)g(X,U)-g(X,Z)g(Y,U)\\+g(JY,Z)g(JX,U)-g(JX,Z)g(JY,U)-2g(JX,Y)g(JZ,U)),\endgather $$
the tensor $\Phi$ is as follows
$$\gather \Phi(X,Y,Z,U)=\frac18(g(Y,Z)h(X,U)-g(X,Z)h(Y,U)\\+g(X,U)h(Y,Z)-g(Y,U)h(X,Z)
+g(JY,Z)\om(X,U)\\-g(JX,Z)\om(Y,U)+g(JX,U)\om(Y,Z)-g(JY,U)\om(X,Z)\\
-2g(JX,Y)\om(Z,U)-2g(JZ,U)\om(X,Y)),\endgather$$ and
$$\Psi(X,Y,Z,U)=-\om(X,Y)\om(Z,U)=-(\om\otimes\om)(X,Y,Z,U).$$

\medskip
{\bf 3. HP metrics on the complex projective spaces $\Bbb{CP}^n$.}
 Let  $\phi:\Bbb {CP}^{n+1}-\{[0,0,...,0,1]\}\r \Bbb{CP}^n$ be a
holomorphic mapping defined as follows:
$$\Bbb {CP}^{n+1}-\{[0,...,0,1]\}\ni[z_0,z_1,...,z_n]\r\phi([z_0,z_1,...,z_n])=   [z_0,z_1,...,z_{n-1}]\in \Bbb{CP}^n.$$
We shall show, that  $\phi$ is the projection on the base of
holomorphic line bundle, whose total space is $H= \Bbb
{CP}^{n+1}-\{[0,0,...,0,1]\}$. Let us consider the mapping
$\phi_i:H_{|U_i} \r U_i\times \Bbb C$
 where  $U_i=\{[z_0,z_1,...,z_{n-1}]\in \Bbb{CP}^n:z_i\ne0\}$, $H_{|U_i}=\phi^{-1}(U_i)$,
 defined in the following way:

 $$\phi_i([z_0,z_1,...,z_n]) = (
 [z_0,z_1,...,z_{n-1}],\frac{z_n}{z_i}).$$
Then
 $ \phi_i\circ \phi^{-1}_j(x , z)=\phi_i([z_0,z_1,...,z_n])$ ,
 where  $x=[z_0,z_1,...,z_{n-1}]$ and $z_n=z_jz$.  Hence
 $\phi_i([z_0,z_1,...,z_n])=(x,\frac{z_n}{z_i})=(x,\frac{z_jz}{z_i})=(x,z\frac{z_j}{z_i})
 $.  It follows that the transition function for our bundle are
 $\phi_{ij}(x)=\frac{z_j}{z_i}$.   It follows that the bundle $H$ is
 isomorphic to the line hyperplane bundle over
 $\Bbb{CP}^n$.  Hence $\Bbb{CP}^{n+1}$ arises from the line hyperplane bundle over
 $\Bbb{CP}^n$ by adding a point to its total space.  We have
 $c_1(\Bbb {CP}^{n-1})=n\a$, where $\a\in H^2(\Bbb {CP}^{n-1},\Bbb
 Z)$ is an indivisible integral class.  Let $p:P\r\Bbb {CP}^{n-1}$ be a circle bundle
over $\Bbb {CP}^{n-1}$ classified by the class $\a$. $P$ is a
principal $S^1$ bundle and let $\th$ be a connection form of $P$.
Then  $[\frac{d\th}{2\pi}]=p^*\a$ in $H^2(\Bbb {CP}^{n-1},\Bbb
R)$. It follows that $\Bbb{CP}^n$ can be described as the product
$[0,L]\times P$ with an equivalence relation, where $p:P\r
\Bbb{CP}^{n-1}$ is the $S^1$ bundle related to the hyperplane
bundle over $\Bbb{CP}^{n-1}$, where $\{0\}\times P$ is identified
to a point $[0,0,...,1]$ and two points $ (L,p), (L,q)$ are
related if $p(p)=p(q)$. Let $g_{\Bbb{CP}^{n-1}}$ be a Fubini-Studi
metric on $\Bbb{CP}^{n-1}$ and $\th$ is the standard connection
form on $P$. A metric
$$g=  dt^2+f(t)^2\th\otimes\th+ h(t)^2 g_{\Bbb{CP}^{n-1}}$$
on the product $[0,L]\times P$ extends to a smooth metric on
$\Bbb{CP}^n$ if $f,h$ are positive on $(0,L)$ which are odd
functions at $0$ satisfying conditions $f(0)=h(0)=0,
f'(0)=h'(0)=1$ and $f$ is an add
 function at $L$ satisfying $f(L)=0, f'(L)=-1$ and $h$ is an even
 function at $L$ such that $h(L)\ne 0$. This metric is K\"ahler if
 $f=hh'$ and
also admits a holomorphic Killing vector field with a
K\"ahler-Ricci potential $h^2$ (see [J-1],[D-M-1],[D-M-2]).
\medskip
 {\bf Theorem.}  {\it Let us consider an
 analytic real
 function $P$ on $\Bbb{R}$, which is positive on $[0,1)$, even at $0$ and such
 that $P(0)=1,P'(0)=0,P(1)=0,P'(1)=-2$.  Let us consider a
 function $h$ satisfying an equation $h'=\sqrt{P(h)}$ and such
 that $h''=\frac12P'(h), h(0)=0, h'(0)=1$.  Then
$$g=  dt^2+(h'(t)h(t))^2\th\otimes\th+ h(t)^2 g_{\Bbb{CP}^{n-1}}$$
extends to a smooth, K\"ahler metric  on  $\Bbb{CP}^n$, which is a
QCH K\"ahler  metric on the submanifold
$\Bbb{CP}^n-\{[0,0,..,1]\}$.}
\medskip
{\it Proof.} We shall show that $h$
 is an odd function at $0$.  It suffices to show that $h^{(
 2k)}(0)=0$ for every $k\in\Bbb{N}$.  For $k=0,1$ this equality holds
 true.    Let us assume that it holds true for $l<k$.
 Note that    $h^{(3)}=\frac12P''(h)h'$.  Consequently
 $$ 2h^{(2k)}(0)=(P''(h)h')^{(2k-3)}(0).$$
We first show that $\frac{d^l}{dt^l}(P(h))(0)=0$ for an odd
$l<2k$. It holds true for $k=1$, since
$\frac{d}{dt}(P(h))(0)=P'(0)h'(0)=0$.  Next
$\frac{d^l}{dt^l}(P(h))(0)=\frac{d^{l-1}}{dt^{l-1}}(P'(h)h')= \sum
C^p_{l-1}\frac{d^{p}}{dt^{p}}(P'(h))(0)(\frac{d^{l-1-p}}{dt^{l-1-p}}h)(0)=0.$
Hence
$$2h^{(2k)}(0)=\sum
C^l_{2k-3}\frac{d^l}{dt^l}(P^{(2)}(h))(0)h^{(2k-2-l)}(0)=0,$$
where $C^l_{2k-3}=\frac{(2k-3)!}{l!(2k-3-l)!}$ since for an odd
$l$ we have $\frac{d^l}{dt^l}(P^{(2)}(h))(0)=0$  and for an even
$l$ we have by an induction assumption $h^{(2k-2-l)}(0)=0$. It
follows that if $P$ is an analytic
 function  which is positive on $[0,1)$, even at $0$ and such
 that $$P(0)=1,P'(0)=0,P(1)=0,P'(1)=-2\tag 3.1$$ and $h$ satisfies an
 equation  $h''=\frac12P'(h), h(0)=0, h'(0)=1$ then
 $h'=\sqrt{P(h)}$ and a metric
$$g=  dt^2+(h'(t)h(t))^2\th\otimes\th+ h(t)^2 g_{\Bbb{CP}^{n-1}}$$
is a K\"ahler metric on $\Bbb{CP}^n$, which is a QCH metric on a
manifold $\Bbb{CP}^n-\{[0,0,..,1]\}$  (see [J-1], note that we
write $h= r\sqrt{n}$ and $s=\frac 2n$ since in our case $k=1$).
Hence this metric is a HP metric on a dense, open subset of
$\Bbb{CP}^n$, which means that is a HP metric on the whole of
$\Bbb{CP}^n$.  In fact $R.R=\phi\Pi.R,$ where
$\phi=-4\frac{h''}h=-\frac{2P'(h)}{h}$ (see [J-2]).  Note that  in
[J-2] there is a sign mistake and the formula for $a+\frac b2$
should be
$$ a+\frac b2=
4(\frac{(r')^2}{r^2}-\frac{f'r'}{fr})=-4\frac{r''}r.$$ However the
fact that $a+\frac b2$ changes sign in the case considered in
[J-2] remains true. The function $\phi=a+\frac b2$ depends only on
$t$ and extends smoothly to the whole of $\Bbb{CP}^n$. We have
$$\phi([0,0,..,0,1])=-\lim_{h\r0}\frac{2P'(h)}{h}=-2P''(0).$$$\k$
\medskip
Let us consider as an example a family of polynomials
$$P_{\a}(t)=1+(\a-1)t^2-2\a t^4+\a t^6.$$ Note that for $\a>-4$ every polynomial $P_{\a}$
is positive on the interval $[0,1)$ and satisfies  conditions
$(3.1)$. Let $\a>-4$ and $h_{\a}$ be a solution of a problem
$$h_{\a}''=\frac12P_{\a}'(h_{\a}), h_{\a}(0)=0, h_{\a}'(0)=1.$$ Then we obtain a
family of $HP$-metrics
$g_{\a}=dt^2+(h_{\a}'(t)h_{\a}(t))^2\th\otimes\th+
h_{\a}(t)^2g_{\Bbb{CP}^{n-1}}$ on $\Bbb{CP}^n$.  Note that for
$\a=0$ we get $h_0(t)=\sin t$ and we obtain a standard symmetric
metric on $\Bbb{CP}^n$ (see [P], p. 17).  We also have
$\phi_{\a}=-\frac{2P_{\a}'(h_a)}{h_{\a}}=4(-(\a-1)+4\a
h_{\a}^2-3\a h_{\a}^4)$. If $\a\in(-3,1)$ then $\phi_{\a}>0$ and
consequently we get examples of compact HP K\"ahler manifolds with
$\phi>0$. If $\a\in\{-3,1\}$ then $\phi_{\a}\ge 0$.  Z.  Olszak in
[O-2] proved that a compact HP K\"ahler manifold which has
$\phi\ge0$ and constant scalar curvature must be locally
symmetric.  Our examples show that the assumption of constant
scalar curvature in Olszak's theorem is necessary.

\bigskip
\centerline{\bf References.}

\medskip
[D] R. Deszcz ,  {\it On pseudosymmetric spaces},  Bull. Soc.
Math. Belg., S\'er.A, (1992),  44, 1-34.
\medskip
\cite{D-M-1} A. Derdzi\'nski, G. Maschler {\it Special
K\"ahler-Ricci potentials on compact K\" ahler manifolds}, J.
reine angew. Math. 593 (2006), 73-116.
\medskip
\cite{D-M-2}  A. Derdzi\'nski, G.  Maschler {\it Local
classification of conformally-Einstein  K\"ahler metrics in higher
dimensions}, Proc. London Math. Soc. (3) 87 (2003), no. 3,
779-819.

\medskip
[H] Hotlo\'s M.,  {\it On holomorphically pseudosymmetric
K\"ahlerian manifolds}, In: Geometry and topology of submanifolds,
VII (Leuven 1994,Brussels 1994), 139-142, World Sci. Publ. River
Edge, NJ, (1995)
\par
\medskip
\cite{G-M-1} G.Ganchev, V. Mihova {\it K\"ahler manifolds of
quasi-constant holomorphic sectional curvatures}, Cent. Eur. J.
Math. 6(1),(2008), 43-75.
\par
\medskip
\cite{G-M-2} G.Ganchev, V. Mihova {\it Warped product K\"ahler
manifolds and Bochner-K\"ahler metrics}, J. Geom. Phys. 58(2008),
803-824.
\par
\medskip
 [J-1]  W.  Jelonek , {\it K\"ahler manifolds with
quasi-constant holomorphic curvature}, Ann. Global Analysis and
Geom. (2009),36,143-159, DOI : 10.1007/s10455-009-9154-z.
\medskip
[J-2]  W.  Jelonek, {\it Compact holomorphically pseudosymmetric
K\"ahler manifolds}, Colloquium Mathematicum, vol.117,No.2, (2009)
243-249.
\medskip
[O-1] Z.  Olszak ,  {\it Bochner flat K\"ahlerian manifolds with a
certain condition on the Ricci tensor}, Simon Stevin, (1989),  63,
295-303.

\medskip
[O-2] Z.  Olszak ,  {\it On compact holomorphically
pseudosymmetric K\"ahlerian  manifolds}, Cent. Eur. J. Math. 7(3),
 2009, 442-451.
\medskip
[P]   P. Petersen , Riemannian Geometry,  Graduate Texts in
Mathematics, Sprin\-ger, 2006
\medskip
[Sz-1] Z.I. Szab\'o,  {\it Structure theorems on Riemannian spaces
satisfying $R(X,Y).R$ $=0$. I. The local version}, J. Diff. Geom.,
(1982), 17, 531-582.

\medskip
[Sz-2] Z.I. Szab\'o,  {\it Structure theorems on Riemannian spaces
satisfying $R(X,Y).R$ $=0$. II. Global versions},  Geom. Dedicata,
(1985), 19, 65-108.
\medskip
[Y]  Yaprak S., {\it Pseudosymmetry type curvature conditions on
K\"ahler hypersurfaces}, Math. J.  Toyama Univ., (1995), 18,
107-136.

\medskip Institute of Mathematics

Technical University of Cracow

 Warszawska 24

31-155 Krak\'ow,POLAND.

E-mail address: wjelon\@pk.edu.pl

\end